\documentclass[12pt,thmsa]{article}
\usepackage{amsfonts}
\begin{document}

\author{David Carf\`{i}}
\title{Summable families in tempered distribution spaces}
\date{}
\maketitle

\begin{abstract}
In this note we define summable families in tempered distribution spaces and
state some their properties and characterizations. Summable families are the
analogous of summable sequences in separable Hilbert spaces, but in tempered
distribution spaces having elements (functional) realizable as generalized
vectors indexed by real Euclidean spaces (not pointwise defined ordered
families of scalars indexed by real Euclidean spaces in the sense of
distributions). Any family we introduce here is summable with respect to
every tempered system of coefficients belonging to a certain normal spaces
of distributions, in the sense of superpositions. The summable families we
present in this note are one possible rigorous and simply manageable
mathematical model for the infinite families of vector-states appearing in
the formulation of the continuous version of the celebrated Principle of
Superpositions in Quantum Mechanics.
\end{abstract}

\bigskip

\section{\textbf{Characterizations of}$\mathcal{\ }^{\mathcal{D}}$\textbf{%
families} (*)}

\bigskip

As we proved the summability theorem of Schwartz families and
characterization of summability, in a perfectly analogous way, it can be
proved the following theorem.

\bigskip

The $^{\mathcal{D}}$families in $\mathcal{D}_{n}^{\prime }$ can be defined
analogously to the $^{\mathcal{S}}$families in $\mathcal{S}_{n}^{\prime }$.

\bigskip

\textbf{Definition}\emph{\ }\textbf{(family of tempered distributions of
class\ }$\mathcal{D}$\textbf{).}\emph{\ Let }$v$\emph{\ be a family of
distributions in the space }$\mathcal{D}_{n}^{\prime }$\emph{\ indexed by
the Euclidean space }$\Bbb{R}^{m}$\emph{. The family }$v$\emph{\ is called a 
\textbf{Schwartz family of class} }$\mathcal{S}$\emph{\ or even }$^{\mathcal{%
D}}$\emph{\textbf{family} if, for each test function }$\phi \in \mathcal{D}%
_{n}$\emph{, the image of the test function }$\phi $\emph{\ by the family }$%
v $\emph{\ - that is the function }$v(\phi ):\Bbb{R}^{m}\rightarrow \Bbb{K}$%
\emph{\ defined by} 
\[
v(\phi )(p):=v_{p}(\phi ), 
\]
\emph{for each index }$p\in \Bbb{R}^{m}$\emph{\ - belongs to the space of
test functions }$\mathcal{D}_{m}.$\emph{\ We shall denote the set of all }$^{%
\mathcal{D}}$\emph{families by }$\mathcal{D}(\Bbb{R}^{m},\mathcal{D}%
_{n}^{\prime })$\emph{.}

\bigskip

The following theorem holds since the Corollary of page 91 of
Dieudonn\'{e}-Schwartz seminal paper holds true because $\mathcal{D}%
_{n}^{\prime }$ is an $\mathcal{LF}$-space.

\bigskip

\textbf{Theorem (basic properties on }$^{\mathcal{D}}$\textbf{families}).%
\emph{\ Let }$v\in \mathcal{D}(\Bbb{R}^{m},\mathcal{D}_{n}^{\prime })$\emph{%
\ be a family of distributions. Then the following assertions hold and are
equivalent:}

\begin{itemize}
\item[\emph{1)}]  \emph{\ for every }$a\in \mathcal{D}_{m}^{\prime }$\emph{\
the composition }$u=a\circ \widehat{v}$\emph{, i.e., the functional} 
\[
u:\mathcal{D}_{n}\rightarrow \Bbb{K}:\phi \mapsto a\left( \widehat{v}(\phi
)\right) , 
\]
\emph{is a distribution;}

\item[\emph{2)}]  \emph{\ the operator }$\widehat{v}$\emph{\ is transposable;%
}

\item[\emph{3)}]  \emph{\ the operator }$\widehat{v}$\emph{\ is }$\left(
\sigma (\mathcal{D}_{n}),\sigma (\mathcal{D}_{m})\right) $\emph{-continuous
from }$\mathcal{D}_{n}$\emph{\ to }$\mathcal{D}_{m}$\emph{;}

\item[\emph{4)}]  \emph{\ the operator }$\widehat{v}$\emph{\ is a strongly
continuous from }$\left( \mathcal{D}_{n}\right) $\emph{\ to }$\left( 
\mathcal{D}_{m}\right) $\emph{.}
\end{itemize}

\bigskip

\section{\textbf{Algebraic }$^{E}$\textbf{Families and }$^{E}$\textbf{%
summable families}}

\bigskip

Let us begin with a family of tempered distribution which is not of class $%
\mathcal{S}$.

\bigskip

\textbf{Example (a family\ that is not of class }$\mathcal{S}$\textbf{).}
Let $u$ be a distribution in $\mathcal{S}_{n}^{\prime }$ and let $v$ be the
constant family in $\mathcal{S}_{n}^{\prime }$, indexed by the Euclidean
space $\Bbb{R}^{m}$, defined by $v_{y}=u$, for each point $y\in \Bbb{R}^{m}$%
. Then, if the distribution $u$ is different from zero, $v$ is not of class $%
\mathcal{S}$. In fact, let $\phi \in \mathcal{S}(\Bbb{R}^{n},\Bbb{K})$ be
such that $u(\phi )\neq 0$; for every point-index $y\in \Bbb{R}^{m}$, we
have 
\begin{eqnarray*}
v(\phi )(y) &=&v_{y}(\phi )= \\
&=&u(\phi )1_{\Bbb{R}^{m}}(y),
\end{eqnarray*}
where, $1_{\Bbb{R}^{m}}$ is the constant $\Bbb{K}$-functional on $\Bbb{R}%
^{m} $ with value $1$. Thus, the function $v(\phi )$ is a constant $\Bbb{K}$%
-functional on $\Bbb{R}^{m}$ different from zero, and so it cannot live in
the space $\mathcal{S}_{m}$.

\bigskip

The preceding example induces us to consider other classes of families in
addition to the $^{\mathcal{S}}$families, for this reason, we will give the
following definitions.

\bigskip

We shall denote by $\mathcal{C}_{m}$ the space $\mathcal{C}^{0}(\Bbb{R}^{m},%
\Bbb{K})$ of continuous functions defined on the Euclidean space $\Bbb{R}%
^{m} $ and with values in the scalar field $\Bbb{K}$.

\bigskip

\textbf{Definition (}$^{E}$\textbf{families and algebraically }$^{E}$\textbf{%
summable families).} \emph{Let }$E$\emph{\ be a subspace of the function
space }$\mathcal{F}(\Bbb{R}^{m},\Bbb{K})$\emph{\ (without any topology)
containing the space }$\mathcal{S}_{m}$\emph{.\ If }$v$\emph{\ is a family
in the distribution space }$\mathcal{S}_{n}^{\prime }$\emph{\ indexed by }$%
\Bbb{R}^{m}$\emph{, we say that the family }$v$\emph{\ is an }$^{E}$\emph{%
\textbf{family} if, for every test function }$\phi $\emph{\ in }$\mathcal{S}%
_{n}$\emph{, the image }$v(\phi )$\emph{\ of the test function by the family 
}$v$\emph{\ lies in the subspace }$E$\emph{. An }$^{E}$\emph{family }$v$%
\emph{\ is said to be \textbf{algebraically} }$^{E}$\emph{\textbf{summable}
or }$^{E^{*}}$\emph{summable if, for every functional }$a$\emph{\ in the
algebraic dual }$E^{*}$\emph{, the functional} 
\[
\int_{\Bbb{R}^{m}}av:\mathcal{S}_{n}\rightarrow \Bbb{K}:\phi \mapsto
a(v(\phi )) 
\]
\emph{is a tempered distribution living in }$\mathcal{S}_{n}^{\prime }$\emph{%
. More generally, if }$F$\emph{\ is a part of the algebraic dual }$E^{*}$%
\emph{\ we say that the family }$v$\emph{\ is }$^{F}$\emph{summable if the
above functional }$\phi \mapsto a(v(\phi ))$ \emph{is a tempered
distribution living in }$\mathcal{S}_{n}^{\prime }$\emph{, for every
functional }$a$\emph{\ in }$F$\emph{.}

\bigskip

\textbf{Example.} With the preceding new definition, the family of the above
example is a $^{\mathcal{E}_{m}}$family in $\mathcal{S}_{n}^{\prime }$,
where by $\mathcal{E}_{m}$ we (in standard way) denote the space $C^{\infty
}(\Bbb{R}^{m},\Bbb{K})$ of smooth function from $\Bbb{R}^{m}$ into $\Bbb{K}$.

Moreover, for every tempered distribution $a$ in $\mathcal{E}_{m}^{\prime }$%
, we have 
\begin{eqnarray*}
a(v(\phi )) &=&a(u(\phi )1_{\Bbb{R}^{m}})= \\
&=&u(\phi )a(1_{\Bbb{R}^{m}})= \\
&=&u(\phi )\int_{\Bbb{R}^{m}}a,
\end{eqnarray*}
where by $\int_{\Bbb{R}^{m}}a$ we denote the integral of the distribution $a$
(we recall that the compact support distributions are integrable and their
integrals is defined as their value on the constant unit functional $1_{\Bbb{%
\ R}^{m}}$), so that 
\[
\int_{\Bbb{R}^{m}}av=\left( \int_{\Bbb{R}^{m}}a\right) u, 
\]
and the family $v$ is $^{\mathcal{E}_{m}^{\prime }}$summable.

\bigskip

\section{$^{E}$\textbf{Families and }$^{E}$\textbf{summable families}}

\bigskip

\textbf{Remark.} Let $E$\ be a subspace of the function space $\mathcal{F}(%
\Bbb{R}^{m},\Bbb{K})$\ (without any topology) containing the space $\mathcal{%
\ S}_{m}$ and let $w$\ be a Hausdorff locally convex topology on the
subspace $E$.

\begin{itemize}
\item  If the topological vector space $\left( \mathcal{S}_{m}\right) $\ is
continuously imbedded in the space $E_{w}$, then, the topological dual $%
E_{w}^{\prime }$\ is continuously imbedded in the space $\mathcal{S}%
_{m}^{\prime }$. In this case, in Distribution Theory, we say that the dual $%
E_{w}^{\prime }$ is \emph{a space of tempered distribution on} $\Bbb{R}^{m}$.

\item  Moreover, if the topological vector space $E_{w}$ is continuously
imbedded in the space $\left( \mathcal{C}_{m}\right) $, then the dual $%
\mathcal{C}_{m}^{\prime }$ is contained in the dual $E_{w}^{\prime }$. In
other terms, every Radon measure on $\Bbb{R}^{m}$ with compact support
belongs to the dual $E_{w}^{\prime }$ and, in particular, the Dirac family
is contained in the topological dual $E_{w}^{\prime }$. Since the Dirac
basis in sequentially total in the space of tempered distributions $\mathcal{%
S}_{m}^{\prime }$ (the linear hull of the Dirac basis is dense in $\mathcal{S%
}_{m}^{\prime }$) and since the space $E_{w}^{\prime }$ is continuously
imbedded in the space $\mathcal{S}_{m}^{\prime }$ itself, the Dirac family
shall be sequentially total also in the topological vector space $%
(E_{w}^{\prime })_{\sigma }$, that is sequentially dense with respect to the
weak* topology $\sigma (E^{\prime },E)$.
\end{itemize}

\bigskip

Now we can give two new definitions.

\bigskip

\textbf{Definition (}$^{E}$\textbf{families and }$^{E}$\textbf{summable
families).} \emph{Let }$E$\emph{\ be a subspace of the space }$\mathcal{F}(%
\Bbb{R}^{m},\Bbb{K})$\emph{\ containing the space }$\mathcal{S}_{m}$\emph{\
and endowed with a locally convex linear topology }$w$\emph{.\ If }$v$\emph{%
\ is a family in the distribution space }$\mathcal{S}_{n}^{\prime }$\emph{\
indexed by }$\Bbb{R}^{m}$\emph{. We say that the family }$v$\emph{\ is an }$%
^{E}$\emph{\textbf{family} if, for every test function }$\phi $\emph{\ in }$%
\mathcal{S}_{n}$\emph{, the image }$v(\phi )$\emph{\ of the test function by
the family }$v$\emph{\ lies in the subspace }$E$\emph{. An }$^{E}$\emph{%
family is said to be }$^{E}$\emph{\textbf{summable} if for every tempered
distribution }$a$\emph{\ in the topological dual }$E_{w}^{\prime }$\emph{\
the functional} 
\[
\mathcal{S}_{n}\rightarrow \Bbb{K}:\phi \mapsto a(v(\phi )) 
\]
\emph{is a tempered distribution in }$\mathcal{S}_{n}^{\prime }$\emph{.}

\bigskip

\section{\textbf{Normal spaces of distributions and summability}}

\bigskip

\textbf{Definition (of normal space of test function for }$\mathcal{S}
_{m}^{\prime }$\textbf{).} \emph{We will call a locally convex topological
vector space }$E$\emph{\ a \textbf{normal space of test functions for the
distribution space }}$\mathcal{S}_{m}^{\prime }$\emph{\ if it verifies the
following properties}

\begin{itemize}
\item  \emph{\ the space }$E$\emph{\ is an algebraic subspace of the space }$%
\mathcal{C}_{m}$\emph{;}

\item  \emph{\ the space }$E$\emph{\ contains the space }$\mathcal{S}_{m}$%
\emph{;}

\item  \emph{\ the topological vector space }$\left( \mathcal{S}_{m}\right) $%
\emph{\ is continuously imbedded and dense in the topological vector space }$%
E$\emph{;}

\item  \emph{\ the topological vector space }$E$\emph{\ is continuously
imbedded in the space }$\left( \mathcal{C}_{m}\right) $\emph{.}
\end{itemize}

\emph{In these conditions the dual }$E^{\prime }$\emph{\ is called \textbf{a
normal space of tempered distributions on} }$\Bbb{R}^{m}$\emph{.}

\bigskip

\textbf{Theorem (on the }$^{E}$\textbf{family generated by a linear and
continuous operator).}\emph{\ Let }$E$\emph{\ be a normal space of test
functions for the space }$\mathcal{S}_{m}^{\prime }$\emph{, let }$A:\mathcal{%
S}_{n}\rightarrow E$\emph{\ be a linear and continuous operator of the space 
}$(\mathcal{S}_{n})$\emph{\ into the space }$E$\emph{\ and let }$\delta $%
\emph{\ be the Dirac family in }$\mathcal{C}_{m}^{\prime }$\emph{. Then, the
family of functionals} 
\[
A^{\vee }:=(\delta _{p}\circ A)_{p\in \Bbb{R}^{m}} 
\]
\emph{is a family of distribution in }$\mathcal{S}_{n}^{\prime }$\emph{\ and
it is an }$^{E}$\emph{family.}

\bigskip

We can prove that:

\bigskip

\textbf{Theorem. }\emph{Let }$E$\emph{\ be a normal space of test functions
for the distribution space }$\mathcal{S}_{m}^{\prime }$\emph{. Then, every }$%
^{E}$\emph{family in }$\mathcal{S}_{n}^{\prime }$\emph{\ (obviously indexed
by the }$m$\emph{-dimensional Euclidean space) is }$^{E}$\emph{summable.}

\bigskip

\textbf{David Carf\`{i}}

\emph{Faculty of Economics}

\emph{University of Messina}

\emph{davidcarfi71@yahoo.it}

\end{document}